\documentclass[12pt,a4paper]{amsart}

\usepackage {amsfonts}
\usepackage[english]{babel}
\def\g{\gamma}
\def\G{\Gamma}
\def\d{\delta}
\def\p{\varphi}
\def\e{\varepsilon}
\def\l{\lambda}
\def\L{\Lambda}
\def\a{\alpha}

\def\s{\sigma}
\def\R{{\mathbb R}}
\def\Z{{\mathbb Z}}
\def\N{{\mathbb N}}
\def\C{{\mathbb C}}

\DeclareMathOperator{\supp}{supp}

\newtheorem{theorem}{Theorem}

\begin{document}

\title{Uniqueness theorem for Fourier transformable measures on LCA groups}

\author{S.Yu.Favorov}

\address{Sergii Favorov,
\newline\hphantom{iii}  Karazin's Kharkiv National University
\newline\hphantom{iii} Svobody sq., 4,
\newline\hphantom{iii} 61022, Kharkiv, Ukraine}
\email{sfavorov@gmail.com}

\maketitle {\small
\begin{quote}
\noindent{\bf Abstract.}
We show that if  points of  supports of two discrete ''not very thick'' Fourier transformable measures on LCA groups tend to one another at infinity and the same is true for the masses at these points,
then these measures coincide. The result is valid for discrete  almost periodic measures on LCA groups too. Also, we show that the result is false for some discrete ''thick'' measures. To do this,
we construct a discrete almost periodic measure on the real axis, whose masses at the points of support tend to zero as these points approach infinity.
\medskip


\end{quote}
}

\section{introduction}
It was proved in \cite{F1} that if  points of  supports of two discrete ''not very thick'' Fourier quasicrystals tend to one another at infinity and the same is true for the masses at these points,
then these quasicrystals coincide. This result is based on the fact that Fourier quasicrystals are almost periodic tempered distributions. Earlier P.Kurasov and R.Suhr \cite{KS} showed  that
if zeros of two holomorphic almost periodic functions in a strip get closer at infinity, then the zeros sets of these functions coincide. It is natural to expect the same effect
for other almost periodic objects, in particular, for discrete transformable measures on LCA groups (see, for example, \cite{LA} or \cite{MS}), which in a sense can be regarded as  analogs
of quasicrystals for LCA groups.
\smallskip

To formulate our results we need some definitions.

Let $G$ be a  locally compact second countable $\sigma$-compact Abelian group, which is not compact, $\hat G$ be the group of characters on $G$, $C_u(G)$ be the space of uniformly continuous bounded functions on $G$ with sup-norm $\|.\|_\infty$,
 $C_c(G)$ be its subspace of all compactly supported continuous functions, and $K_2(G)$ be the span of the set $\{f_1\star f_2:\  f_1, f_2\in C_c(G)\}$.
 We will consider Radon measures on $G$, i.e., complex-valued $\sigma$-additive measures $\mu$ on the algebra of Borel sets in $G$ whose variation $|\mu|$ is bounded on every compact subset of $G$
and  for any $\e>0$ and  Borel set $A\subset G$ there is a compact set $K_\e\subset A$ such that  $|\mu|(A\setminus K_\e)<\e$.
We will say that  a measure $\mu$ on $G$ with support $\L$ is {\it discrete}, if for each $\l\in\L$ there is a  neighborhood $\l+U$ of $\l$ such that  $(\l+U)\cap\L=\{\l\}$.
If this is the case, we will write
$$
\mu=\sum_{\l\in\L}\mu(\l)\d_\l,\quad \mu(\l)\in\C,
$$
where $\d_\l$ means the unit mass at the point $\l$. A measure $\mu$ on $G$ is called translation bounded if for all compact sets $K\subset G$ we have
 \begin{equation}\label{tr}
 \sup_{t\in G}|\mu|(K+t)<\infty.
 \end{equation}
 A measure $\mu$ is {\it Fourier transformable} if there is a measure $\hat\mu$ on $\hat G$ such that for all $g\in K_2(G)$
$$
i)\quad (\mu,g)=(\hat\mu,\check g),
$$
$$
ii)\quad \hat g\in L^1(|\hat\mu|).
$$
\noindent Here $\hat g$ is the Fourier transform of $g$, and $\check g$ is the inverse Fourier transform.

Fourier transformable discrete  measures were considered in a series of papers \cite{MS}-\cite{S4}.

A function $g\in C_u(G)$ is {\it almost periodic} if the closure of the family of translates $\{g(\cdot-t)\}_{t\in G}$ is a compact subset of $C_u(G)$. A Radon measure $\mu$ on $G$ is {\it almost periodic}
 if for each function $f\in C_c(G)$ the function $f\star\mu$ is almost periodic on $G$.

By $\# A$ denote a number of points of the finite set $A$.

Denote by $DM$ the class of discrete Radon measures $\mu$ on $G$ with the property
\begin{equation}\label{spar}
\sup_{x\in G}\#\{(x+U)\cap\supp\mu\}<\infty
\end{equation}
for some neighborhood $U$ of $0$.
Also, we  say that discrete measures
$$
\mu=\sum_{\l\in\L} \mu(\l)\d_{\l}\quad\mbox{ and }\quad \nu=\sum_{\g\in\G} \nu(\g)\d_{\g}
$$
 come close to infinity, if there is a bijection $\s$
between  $\L=\supp\mu$ and $\G=\supp\nu$  with the property: for any $\e>0$ and any neighborhood $V$ of $0$ there is a compact $K$ such that
$$
  \l\in\L\setminus K \mbox{ or } \s(\l)\in\G\setminus K\ \Longrightarrow\ \l-\s(\l)\in V \mbox{ and } |\mu(\l)-\nu(\s(\l))|<\e.
$$
\begin{theorem}\label{T1}
If translation bounded Fourier transformable measures $\mu,\ \nu\in DM$ with supports $\L,\G$, respectively, have pure point measures $\hat\mu, \hat\nu$ and come close to infinity,
then the measures $\mu,\,\nu$  coincide.
\end{theorem}

Note that translation bounded Fourier transformable measure $\mu$ is almost periodic if and only if $\hat\mu$ is pure point measure (see \cite{LA}, \cite{MS}). Therefore, the above theorem
 follows from the following one

\begin{theorem}\label{T2}
If almost periodic translation bounded measures $\mu,\ \nu\in DM$  with supports $\L,\G$, respectively, come close to infinity, then they  coincide.
\end{theorem}
Points $\l\in\L,\ \g\in\G$ can be combined into lumps:
\begin{theorem}\label{T3}
Let supports $\L,\,\G$ of discrete almost periodic translation bounded measures $\mu,\ \nu$  decompose into lumps
$$
\L=\cup_{\a\in A}\L_\a,\quad \G=\cup_{\a\in A}\G_\a,\quad  \L_\a\cap\L_{\a'}=\emptyset, \quad \G_\a\cap\G_{\a'}=\emptyset\quad \forall \a\neq\a',
$$
such that for some neighborhood $U$ of $0$
\begin{equation}\label{spar1}
\sup_{x\in G}\#\{\a:\,(x+U)\cap(\L_a\cup\G_\a)\}<\infty,
\end{equation}
and  for every neighborhood $V$  and $\e>0$ there is a compact set $K\subset G$ such that the implication is valid
\begin{equation}\label{c1}
  (\L_a\cup\G_\a)\cap K=\emptyset\ \Longrightarrow\ (\L_a\cup\G_\a)\subset x+V\quad\forall x\in(\L_a\cup\G_\a)\ \mbox{ and}\quad|\mu(\L_\a)-\nu(\G_\a)|<\e,
\end{equation}
then  the measures $\mu$ and $\nu$ coincide.
\end{theorem}
Next, we show that the theorems may be wrong for  discrete measures outside of $DN$.
\begin{theorem}\label{T4}
There is a positive discrete translation bounded almost periodic measure $\mu$ on $G=\R$ such that
\begin{equation}\label{l}
\mu(\l)\to 0\quad\mbox{ as }\quad \l\to\infty,\ \l\in\supp\mu.
\end{equation}
\end{theorem}
We see that the discrete measures $\mu$ and $2\mu$ are almost periodic, translation bounded, and come close to infinity, but do not coincide.

\section{Proof of Theorem \ref{T2}}

 Assume that the measures $\mu$ and $\nu$ do not coincide. Then there is a point $a\in G$ such that $\mu(a)\neq\nu(a)$. Without loss of generality suppose that $a=0$.
Take a compact symmetric neighborhood $U$ of zero such that
$$
U\cap(\L\cup\G)=\{0\}.
$$
Using (\ref{spar}) and reducing $U$, if necessary, we find $N\in\N$ such that for all $x\in G$
\begin{equation}\label{N}
  \#\{(x+U)\cap\L\}<N,\qquad \#\{(x+U)\cap\G\}<N.
\end{equation}
Furthermore, since the function
$$
x_1+x_2+\dots +x_{3N+2},\quad x_j\in G\ \ \forall j,
$$
 is continuous on $G^{3N+2}$, and $0\in G$ has a basis of compact and symmetric neighborhoods,
we get that there exists a compact neighborhood $V$ of zero such that $V=-V$ and all algebraic sums
$$
V_k=V+V+\dots+V\quad (k\mbox{ times},\quad k=1,\dots, 3N+2)
$$
 are compact and symmetric subsets of $U$.
Put
$$
 \p_j(x)=m(V)^{-1}(p_{V}\star p_{V_{3j}})(x),\quad j=1,\dots,N.
$$
Here $m(V)$ means the Haar measure of the set $V$, and
$$
p_A(x)=1 \quad x\in A, \qquad\quad p_A(x)=0 \quad x\not\in A.
$$
Clearly, for all $j=1,\dots,N$ we have $\p_j\in C_c(G)$ and
\begin{equation}\label{p}
  \p_j(x)=1 \mbox{ for } x\in V_{3j-1},\quad \p_j(x)=0 \mbox{ for } x\not\in V_{3j+1},\quad 0\le \p_j(x)\le1 \quad \forall\  x\in G.
\end{equation}
The functions $(\mu-\nu)\star\p_j$ are almost periodic and  $(\mu\star\p_j-\nu\star\p_j)(0)=\mu(0)-\nu(0)\neq0$ for all $j$.
Also, we have
$$
 |(\mu-\nu)\star\p_j(x)|=\left|\sum_{\l\in x+V_{3j+1}}\p_j(x-\l)\mu(\l)-\sum_{\g\in x+V_{3j+1}}\p_j(x-\g)\nu(\g)\right|\le|\mu|(x+U)+|\nu|(x+U).
$$
Since $U$ is a compact subset of $G$,  it follows from (\ref{tr})  that for some $C<\infty$ and for all $j$
 \begin{equation}\label{b}
 \sup_{x\in G}|((\mu-\nu)\star\p_j)(x)|\le C.
 \end{equation}
Put
$$
\Psi(x)=\prod_{j=1}^N ((\mu-\nu)\star\p_j)(x).
$$
This function is almost periodic on $G$ and
\begin{equation}\label{0}
\Psi(0)= (\mu(0)-\nu(0))^N\neq0.
\end{equation}
 Fix $b\in G$, suppose that $\l-\g\in V$ and for some $i,\ 1\le i\le N$,
 $$
 \p_i(b-\l)\neq\p_i(b-\g).
 $$
 Then we get either
 $$
 \l\in b+(V_{3i+1}\setminus V_{3i-1}),\qquad\g\in b+(V_{3i+2}\setminus V_{3i-2}),
 $$
 or
 $$
 \g\in b+(V_{3i+1}\setminus V_{3i-1}),\qquad\l\in b+(V_{3i+2}\setminus V_{3i-2}).
 $$
  Hence, in both cases we have for $j>i$
  $$
  \l,\,\g\in b+V_{3j-1}\quad\mbox{and}\quad \p_j(\l-b)=\p_j(\g-b)=1
  $$
  and for $j<i$
 $$
  \l,\,\g\not\in b+V_{3j+1}\quad\mbox{and}\quad \p_j(\l-b)=\p_j(\g-b)=0.
 $$
  Thus, for every pair $\{\l,\,\s(\l)\}$ such that $\l-\s(\l)\in V$ there is no more than one index $i,\ 1\le i\le N$, such that
 $$
 \p_i(\l-b)\neq\p_i(\s(\l)-b).
 $$
 Also, if $\l$ or $\g$ belong to $b+V_{3N+1}$, then both $\l$ and $\g$ belong to $b+V_{3N+2}\subset U$.
 By (\ref{N}), a number of pairs $\{\l,\,\s(\l)\}$ in the set $U$ is less than $N$, hence there is at least one index $j$ such that
 for every pair
 $$
 \{\l,\,\s(\l)\}\subset b+U,\qquad\l-\s(\l)\in V,\quad\l\in\L,
 $$
we have
$$
 \p_j(b-\l)=\p_j(b-\s(\l)).
$$
 If, besides this, $|\mu(\l)-\nu(\s(\l))|<\e$, then, by (\ref{N}),
$$
 |(\mu\star\p_j-\nu\star\p_j)(b)|=\left|\sum_{\l\in x+V_{3j+1}}\p_j(b-\l)\mu(\l)-\p_j(b-\s(\l))\nu(\s(\l))\right|\le N\e,
$$
and by (\ref{b}),
\begin{equation}\label{P0}
|\Psi(b)|< N\e C^{N-1}.
\end{equation}
Since $\mu$ and $\nu$ come close to infinity, we see that for the neighborhood $V$ there is a compact set $K\subset G$ such that
$$
  \l\in\L\setminus K \mbox{ or } \s(\l)\in\G\setminus K\ \Rightarrow\ \l-\s(\l)\in V \mbox{ and}\quad |\mu(\l)-\nu(\s(\l))|<\e.
 $$
Hence, if $b\not\in K+U$, then under conditions
$$
  \l\in\L, \quad \{\l,\s(\l)\}\subset (b+V_{3N+1}),
$$
we have
$$
\l\not\in K,\qquad \l-\s(\l)\in V,\quad\quad |\mu(\l)-\nu(\s(\l))|<\e,
$$
hence we obtain (\ref{P0}).

 Next, let $G=\cup_n K_n$, where $K_n\subset K_{n+1}$ are compact sets. Fix $x\in G$, and suppose $-x\in K_j$. Take a sequence $b_n\not\in K_n+K+U$. We have $b_n+x\not\in K+U$ for all $n\ge j$, therefore, by (\ref{P0}),
 for every $x\in G$ there is $j=j(x)$ such that
\begin{equation}\label{Ps}
|\Psi(x+b_n)|< N\e C^{N-1},\quad n\ge j.
\end{equation}
 By the definition of almost periodicity, the close of the set $\{\Psi(\cdot+b_n)\}_{n\in\N}$ is a compact subset of $C_u(G)$. Therefore, there is $g\in C_u(G)$ such that an infinite number of functions $\Psi(\cdot+b_n)$
 satisfy the condition
 $$
   \sup_{x\in G} |\Psi(x+b_n)-g(x)|<\e.
 $$
 Using (\ref{Ps}), we obtain
 $$
 |g(x)|\le \e(1+NC^{N-1})\quad \forall x\in G.
 $$
  In particular, $|g(-b_n)|\le \e(1+NC^{N-1})$ and $|\Psi(0)|\le \e(2+NC^{N-1})$, which contradicts to (\ref{0}) for sufficiently small $\e$. The theorem is proved.

              \section{Proof of Theorem \ref{T3}}

 By (\ref{spar1}), take $N\in\N$ and a compact symmetric neighborhood $U$ of zero such that
\begin{equation}\label{N1}
N>\sup_{x\in G}\#\{\a:\,(x+U)\cap(\L_\a\cup\G_\a)\neq\emptyset\}.
\end{equation}
As in the previous proof,  suppose that $\mu(0)-\nu(0)\neq0$ and
$$
U\cap(\L\cup\G)=\{0\}.
$$
As above, take  the sets $V, V_k$ and the functions $\p_j(x),\ \Psi(x)$ satisfying  conditions (\ref{p}), (\ref{b}), and (\ref{0}).
Fix $b\in G$,  suppose that
\begin{equation}\label{sm}
\L_\a\cup\G_\a\subset V+x\qquad\forall\ x\in \L_\a\cup\G_\a,
 \end{equation}
 and for some $\a\in A$ and $i,\  1\le i\le N$,
 \begin{equation}\label{i}
 (\L_\a\cup\G_\a)\cap(b+V_{3i+1}\setminus V_{3i-1})\neq\emptyset.
 \end{equation}
 Then we get
 $$
 \L_\a\cup\G_a\subset b+(V_{3i+2}\setminus V_{3i-2}).
 $$
  Hence, we have for $j>i$
  $$
  \L_\a\cup\G_a\subset V_{3j-1}+b\quad\mbox{and }\quad \p_j(b-y)=1\quad\forall y\in\L_\a\cup\G_\a,
  $$
  and for $j<i$
 $$
  (\L_\a\cup\G_\a)\cap (V_{3j+1}+b)=\emptyset \quad\mbox{and }\quad \p_j(b-y)=0\quad\forall\ y\in\L_\a\cup\G_\a.
 $$
 Thus, for every  $\a$ such that (\ref{sm}) satisfies, there is no more than one index $i,\ 1\le i\le N$ such that (\ref{i}) take place.
  Also, if $(\L_\a\cup\G_\a)\cap(b+V_{3N+1})\neq\emptyset$, then
 $$
  \L_\a\cup\G_a\subset b+V_{3N+2}\subset b+U.
 $$
 By (\ref{N1}),
 $$
  \#\{\a:\,(b+V_{3N+2})\cap(\L_a\cup\G_\a)\neq\emptyset\}<N,
 $$
 hence there is at least one index $j$ such that for all $\a\in A$ under condition (\ref{sm}) we have
 $$
 (\L_\a\cup\G_\a)\cap(b+V_{3j+1}\setminus V_{3j-1})=\emptyset\quad\mbox{ and }\quad\forall\l\in\L_a,\ \forall\g\in\G_\a\quad\  \p_j(b-\l)=\p_j(b-\g).
 $$
  If, besides this, $|\mu(\L_\a)-\nu(\G_\a)|<\e$, then, by (\ref{N1}),
$$
 |(\mu\star\p_j-\nu\star\p_j)(b)|=\left|\sum_{\l\in\L\cap(b+V_{3j+1})}\p_j(b-\l)\mu(\l)-\sum_{\g\in \G\cap(b+V_{3j+1})}\p_j(b-\g)\nu(\g)\right|\le N\e,
$$
and, as above,
\begin{equation}\label{P1}
|\Psi(b)|< N\e C^{N-1}.
\end{equation}
Next, take a compact set $K\subset G$ such that implication (\ref{c1}) satisfies  for the same $\e$ and the neighborhood $V$.
Hence, if
$$
b\not\in K+U,\qquad\L_\a\cup\G_\a\subset b+U,
$$
then
$$
(\L_\a\cup\G_\a)\cap K=\emptyset,\quad\L_\a\cup\G_\a\subset V+x\quad\forall\ x\in (\L_\a\cup\G_\a), \quad |\mu(\L_\a)-\nu(\G_\a)|<\e,
$$
and we get (\ref{P1}).
The end of the proof is the same as in the proof of Theorem \ref{T2}.

\section{Proof of Theorem \ref{T4}}

Take $r_k=2^{-(k+1)^2}$ and put
  $$
    T_k=\frac{1}{2k}\left(\sum_{j=1}^k S_{-r_k j/k}+\sum_{j=1}^k S_{r_k j/k }\right),\quad k=1,2,\dots.
  $$
where $S_t$ means the shift by $t$, i.e., for a measure $\mu$ on $\R$ we have $S_t\mu(A)=\mu(A-t)$ for $A\subset\R$.
 Put
 $$
   \mu_0=\d_0,\quad \mu_1=\mu_0+(S_{-1}+S_1)T_1\mu_0=\d_0+(1/2)(\d_{-1-r_1}+\d_{-1+r_1}+\d_{1-r_1}+\d_{1+r_1}),
  $$
  $$
   \quad \mu_k=\mu_{k-1}+(S_{-3^{k-1}}+S_{3^{k-1}})T_k\mu_{k-1},\quad k=2,3,\dots.
 $$
We need to show that a weak limit $\mu$ of the measures $\mu_k$ exists and has all the properties declared in the theorem.

First, for $N\ge 2$ we have
\begin{equation}\label{est}
  \sum_{k=N}^\infty r_k=2^{-N^2}\sum_{i=1}^\infty 2^{N^2-(N+i)^2}<2^{-N^2}\sum_{i=1}^\infty 2^{-2Ni}=2^{-N^2}\frac{1}{2^{2N}-1}<\frac{r_{N-1}}{3(N-1)}.
\end{equation}
  Also note that for any $x,\ x'\in\R$ we get
 $$
\forall y,\,y'\in\supp T_k(\d_x+\d_{x'}):\,y\neq y'\quad\Longrightarrow\quad|y-y'|\ge\min\{|x-x'|-2r_k,\ r_k/k\}.
$$
  Since
  $$
  3\sum_{i=j}^h r_{k_i}<r_s/s\quad\mbox{for}\quad s<k_j,
  $$
  the supports of the measures $\mu_s$ and  $T_{k_h}T_{k_{h-1}}\dots T_{k_{j+1}}T_{k_j} \mu_s$  are  disjoint.  Then inequality (\ref{est}) for $N=1$ implies that $\sum_{k=1}^\infty r_k<1/3$, therefore by induction on $s$,
$$
  \supp\mu_s\subset I_s=\left(\frac{1-3^s}{2}-\frac{1}{3},\,\frac{3^s-1}{2}+\frac{1}{3}\right),
$$
and the restrictions of measures $\mu_j$ to ${I_s}$ for all $j>s$ coincide with $\mu_s$. Hence there exists a weak limit $\mu$ of the measures $\mu_s$. Note that we have
\begin{equation}\label{mu}
\forall\ n\in\Z\quad\sum\nolimits_{|y-n|<1/3}\mu(y)=1\quad\mbox{and}\quad \supp\mu\subset\bigcup\nolimits_{n\in\Z}(n-1/3,n+1/3).
\end{equation}
Besides, if $x\in\supp\mu\setminus I_s$, then $\mu(x)$ is the mass at the point $x$ of some shift of the measure $T_{k_h}T_{k_{h-1}}\cdots T_{k_{j+1}}T_{k_j}\mu_s$,
where $k_h,\,k_{h-1},\dots k_{j+1},\,k_j$ are some integers such that
$$
   k_h>k_{h-1}>\dots>k_{j+1}>k_j>s.
 $$
 Consequently,
$$
\mu(x)\le 1/(2k_j)<1/(2s),
$$
and  we obtain  (\ref{l}).

Furthermore, take a function $f\in C_c(\R)$. Without loss of generality suppose that $\supp f\subset [-1/6,1/6]$, hence the support of the function $f(x-\cdot)$ for every fixed
$x\in\R$ intersects no more than one of the intervals
$(n-1/3,\ n+1/3)$ for $n\in\Z$.  The above definition of almost periodicity for $g\in C_u(\R)$ is equivalent to the usual one (see, for example, \cite{LZ}):
$$
\forall \e>0\quad \exists L=L(\e)\quad \forall a\in\R\quad \exists \tau\in(a,\,a+L)\quad\forall x\in\R\quad |g(x+\tau)-g(x)|<\e.
$$
Therefore  we will check that for any $\e>0$ there is $s\in\N$ such that for all $x\in\R$ and $\tau=p3^s,\ p\in\Z$,
$$
   |(f\star\mu)(x+\tau)-(f\star\mu)(x)|<\e.
$$
 If
 $$
 [x-1/6,\,x+1/6]\cap\bigcup\nolimits_{n\in\Z}(n-1/3,\ n+1/3)=\emptyset,
 $$
 then for any integer $\tau$
 $$
 [x+\tau-1/6,\,x+\tau+1/6]\cap\bigcup\nolimits_{n\in\Z}(n-1/3,\ n+1/3)=\emptyset
 $$
 as well, and
 $$
   (f\star\mu)(x+\tau)=(f\star\mu)(x)=0.
 $$
 Otherwise, there is a unique $n\in\Z$ such that
 $$
 [x-1/6,\,x+1/6]\cap(n-1/3,\ n+1/3)\neq\emptyset.
 $$
If this is the case, then it is enough for all $x\in(n-1/2,n+1/2)$ to check  the estimate
  \begin{equation}\label{e}
    \left|\int_{|z-n|<1/3}f(x-z)(S_\tau\mu)(dz)-\int_{|z-n|<1/3}f(x-z)\mu(dz)\right|<\e.
 \end{equation}
 Obviously, it is enough to prove this inequality  for $n\in I_s$ (with $\e/2$ instead of $\e$).  In this case the restriction of $\mu$ to $(n-1/3,\ n+1/3)$  coincides
 with the restriction to this interval of the measure $\mu_s$. If $n+\tau\in I_N$, then the restriction of $\mu$ to $(n+\tau-1/3,\ n+\tau+1/3)$  coincides
 with the restriction to this interval of the measure
$$
   S_\tau T_{k_h}T_{k_{h-1}}\cdots T_{k_{j+1}}T_{k_j}\mu_s
 $$
for some
$$
  k_j<k_{j+1}<\dots<k_{h-1}<k_h
 $$
 such that $s<k_j,\  k_h\le N$. Set
$$
q=(2k_j)(2k_{j+1})\dots(2k_{h-1})(2k_h),\qquad \eta=r_{k_j}+r_{k_{j+1}}+\dots+r_{k_{h-1}}+r_{k_h},
$$
and for each  point $y\in\supp\mu_s\cap(n-1/3,\,n+1/3)$
$$
A(y)=\{z\in\R:\,|z-y|\le\eta,\ z+\tau\in\supp\mu\}.
$$
Since supports of the measures
$$
\mu_s,\quad T_{k_j}\mu_s,\quad T_{k_{j+1}}T_{k_j}\mu_s,\quad\dots\quad T_{k_{h-1}}\dots T_{k_{j+1}}T_{k_j}\mu_s,\quad T_{k_h}T_{k_{h-1}}\dots T_{k_{j+1}}T_{k_j}\mu_s
$$
are  mutually disjoint, it follows that  $\# A(y)=q$ and
$$
T_{k_h}T_{k_{h-1}}\cdots T_{k_{j+1}}T_{k_j}\d_y=\frac{1}{q}\sum_{z\in A(y)}\d_z.
$$
Therefore,
$$
 \sum_{|z-n|<1/3}f(x-z)\mu(z+\tau)=\sum_{|y-n|<1/3}\frac{\mu(y)}{q}\sum_{z\in A(y)}f(x-z),
$$
and we have
$$
  \sum_{|z-n|<1/3}f(x-z)\mu(z+\tau)-\sum_{|y-n|<1/3}f(x-y)\mu(y)= \sum_{|y-n|<1/3}\frac{\mu(y)}{q}\sum_{z\in A(y)}[f(x-z)-f(x-y)].
$$
By (\ref{est}),  $|z-y|\le\eta<r_s/s$, hence for $s$ sufficiently large $|f(x-z)-f(x-y)|<\e$.
 Using (\ref{mu}), we obtain (\ref{e}). The theorem is proved.

\end{document}